\documentclass[12pt,thmsa]{article}
\usepackage{amsfonts}
\usepackage{amssymb}
\newtheorem{teo}{Theorem}[section]

\newtheorem{cor}[teo]{Corollary}

\newtheorem{obs2}[teo]{Remark}

\newtheorem{tea}{Theorem}[subsection]

\newtheorem{no2}[teo]{Note}

\newtheorem{no3}[tea]{Note}

\newcommand{\Gal}{{\rm Gal}}
\newcommand{\Frob}{{\rm Frob }}
\newcommand{\trace}{{\rm trace}}
\newcommand{\F}{{\mathbb{F} }}
\newcommand{\Q}{{\mathbb{Q} }}
\newcommand{\mod}{{\rm mod}}

\title{From potential modularity to modularity for
 integral Galois representations and rigid Calabi-Yau threefolds}
\author{Luis V. Dieulefait\\ Universitat de Barcelona
  \thanks{
e-mail: LDieulefait@ub.edu }}

\begin{document}

\maketitle
%\tableofcontents
%\vspace{6mm}
%\newpage
%\noindent Running head:
\begin{abstract}
 In a previous article, we have proved a result asserting the existence
 of a compatible family of Galois representations containing a given
 crystalline irreducible odd two-dimensional representation. We apply this
 result, combined with the potential modularity results of Taylor,
  to prove modularity for
 any irreducible crystalline $\ell$-adic odd $2$-dimensional Galois
 representation (with finite ramification set) unramified at $3$
 verifying an ``ordinarity at $3$" easy to check condition,
 with Hodge-Tate weights $\{0, w \}$ such that $2 w <  \ell$ (and $\ell > 3$)
  and such that the traces $a_p$
 of the images of Frobenii verify $\Q(\{ a_p \}) = \Q $. This result applies in particular
 to any
 motivic compatible family of odd two-dimensional Galois representations of $\Gal(\bar{\Q}/\Q)$
 if the motive has rational coefficients, good reduction at $3$, and the ``ordinarity at $3$" condition
  is satisfied.
As a corollary, this proves that all rigid Calabi-Yau threefolds defined over $\Q$
having good reduction at $3$ and satisfying $ 3 \nmid a_3$ are modular.

%Mathematics Subject Classification: 11F80, 11F46
\end{abstract}
%\newpage

\section{The result and its proof}

The main tools in our proof are the existence of families proved in [D], together
with the potential modularity results proved in [T1], [T2]. We will also need a
result from [W] which controls ordinarity for Hilbert modular forms.\\
We will call field of coefficients of a Galois representation the field generated by
the traces of the images of Frobenii elements.

 \begin{teo}
  \label{teo:Taylor} Let $\ell > 3$ be a prime. Let
  $\sigma_{\ell}$ be a two dimensional odd irreducible $\ell$-adic Galois
  representation (of the absolute Galois group of $\Q$, continuous)
  ramified only at $\ell$ and at a finite set of primes $S$ not containing $3$, with field of
  coefficients $\Q$. Assume that
  $\sigma_\ell$ is crystalline at $\ell$, with Hodge-Tate weights  $\{0, w \}$ ($w$ odd).
  Assume also that $\ell > 2w $. Then, if $ a_3 := \trace ( \sigma_\ell (\Frob \; 3) )$ is
  not divisible by $3$, the representation $\sigma_\ell$ is modular.\\
  \end{teo}

  \begin{cor} If an odd, two-dimensional, compatible family of Galois representations attached
   to a motive defined over $\Q$, having rational coefficients and good reduction at
   $3$,  verifies $3 \nmid a_3 := \trace ( \sigma_\ell (\Frob \; 3) )$, $\ell \neq 3$,
    then the family (thus, the motive) is modular. In
   particular, any rigid Calabi-Yau threefold defined over $\Q$ with good reduction
   at $3$ and $3 \nmid a_3$ is modular.
  \end{cor}

  Remark: This modularity criterion for rigid Calabi-Yau threefolds is different
  than
  those obtained in [DM]. In particular, the criteria in [DM] required good
  reduction at $5$ or $7$.\\

  Proof of corollary: Just apply the theorem to the $\ell$-adic representation in the family
   for a sufficiently large prime $\ell$ where the motive has good reduction.\\

   Proof of theorem:

   First, recall that from the existence of a family result in [D], we can insert $\sigma_\ell$
    in a (strongly) compatible family $\{ \sigma_q \}$, which has rational coefficients and
    is unramified at $3$ for any $q \neq 3$.
     Also, from the results in [T2], we know that
     this family is potentially modular, i.e., that when restricted to some totally real number
      field $F$ all representations in the family agree with those attached to some Hilbert modular
       form $h$ over $F$.\\
    Now, we will
   translate the ``easy to check condition" on the trace of $\sigma_\ell$ at
$\Frob \; 3$ into ordinarity of the $3$-adic representation in the compatible
family. We apply a  result of Wiles [W], which tells us that we can read in the
corresponding eigenvalue of a Hilbert modular form that the attached Galois
representation is ordinary. Our family of representations becomes modular when
restricted to the Galois group of a totally real number field F which can be assumed
(using solvable base change) to be such that $3$ is totally split in $F/ \Q$ (cf.
[D]), so that ordinarity at $3$ of the restriction to the Galois group of $F$ is
equivalent to ordinarity of the full $3$-adic representation. The condition $3 \nmid
a_3$ implies (cf. [W]) that when we restrict $\sigma_3$ to the Galois group of $F$
we get a  modular Galois representation which is ordinary at (every prime
of $F$ dividing) $3$, therefore we conclude that $\sigma_3$ is ordinary. \\

 The proof finishes using
Skinner-Wiles results: since the residual $\mod \; 3$ representation has
coefficients in $\F_3$, it is known that it is either modular or reducible (by
results of Langlands and Tunnell). Knowing that $\sigma_3$ is ordinary, an
application of [SW1] and [SW2] gives the modularity of $\sigma_3$, thus of
$\sigma_\ell$ (because they are compatible).\\

\section{Bibliography}

[D] Dieulefait, L., {\it Existence of compatible families and new cases of the
Fontaine-Mazur conjecture}, to appear in J. Reine Angew. Math.; available at
http://front.math.ucdavis.edu/math.NT \\

[DM] Dieulefait, L., Manoharmayum, J., {\it Modularity of rigid Calabi-Yau
threefolds over $\Q$}, in ``Calabi-Yau Varieties and Mirror Symmetry", Fields
Institute Communications, {\bf 38}, AMS (2003)\\

[SW1] Skinner, C., Wiles, A., {\it Residually reducible representations and modular
forms},  Publ. Math. IHES {\bf 89} (2000)\\

[SW2] Skinner, C., Wiles, A., {\it Nearly ordinary deformations of irreducible
residual representations}, Ann. Fac. Sci. Toulouse Math. (6) {\bf 10} (2001)
\\

 [T1] Taylor, R., {\it Remarks on a conjecture
of Fontaine and Mazur}, J. Inst. Math.
Jussieu {\bf 1} (2002)\\

 [T2] Taylor, R., {\it On the meromorphic continuation of
degree two
 L-functions}, preprint, (2001);
 available at http://abel.math.harvard.edu/$\sim$rtaylor/ \\

[W] Wiles, A., {\it On ordinary $\lambda$-adic representations associated to modular
forms},  Invent. Math. {\bf 94} (1988)\\

\end{document}